\documentclass[reqno]{amsart}
\usepackage{hyperref}

\begin{document}

\title[\hfilneg \hfil Uniqueness of  meromorphic function sharing three small functions CM with its $n-$ exact difference]
{Uniqueness of meromorphic function sharing three small functions CM with its $n-$ exact difference}

\author[XiaoHuang Huang \hfil \hfilneg]
{XiaoHuang Huang}

\address{XiaoHuang Huang: Corresponding author\newline
Department of Mathematics, Shenzhen University, Shenzhen 518055, China}
\email{1838394005@qq.com}

\subjclass[2010]{30D35}
\keywords{ Uniqueness, meromorphic functions,  share small functions, differences}
\begin{abstract}
In this paper, we study the uniqueness of the difference of meromorphic functions. We prove the following result: Let $f$ be a  non-constant meromorphic function of hyper-order less than $1$, let $\eta$ be a non-zero complex number, $n\geq1$, an integer, and let $a,b,c\in\hat{S}(f)$  be three  distinct small functions and two of them  be periodic small functions with period $\eta$. If $f$ and $\Delta_{\eta}^{n}f$ share $a,b,c$ CM, then $f\equiv\Delta_{\eta}^{n}f$.
\end{abstract}

\maketitle
\numberwithin{equation}{section}
\newtheorem{theorem}{Theorem}[section]
\newtheorem{lemma}[theorem]{Lemma}
\newtheorem{remark}[theorem]{Remark}
\newtheorem{corollary}[theorem]{Corollary}
\newtheorem{example}[theorem]{Example}
\newtheorem{problem}[theorem]{Problem}
\allowdisplaybreaks

\section{Introduction and main results}
In this paper, we use the standard denotations in the  Nevanlinna value distribution theory, see(\cite{h3,y1,y2}).
 Throughout this paper,  $f(z)$ is a meromorphic function on the whole complex plane.

 By $S(r,f)$, we denote any quantity satisfying $S(r, f) = o(T(r, f))$, as $r\to \infty $ outside of a possible exceptional set of finite logarithmic measure. A meromorphic function $a(z)$ satisfying $T(r,a)=S(r,f)$ is called a small function of $f$.  We denote $S(f)$ as the family of all small meromorphic functions of $f$ which includes the constants in $\mathbb{C}$. Moreover, we define $\hat{S}(f)=S(f)\cup\{\infty\}$. We say that two non-constant meromorphic functions $f$ and $g$ share small function $a$ CM(IM) if $f-a$ and $g-a$ have the same zeros counting multiplicities (ignoring multiplicities).

   Define
 $$\lambda(f)=\varlimsup_{r\rightarrow\infty}\frac{log^{+}T(r,f)}{logr},$$
 $$\rho(f)=\varliminf_{r\rightarrow\infty}\frac{log^{+}T(r,f)}{logr},$$
 $$\rho_{2}(f)=\varlimsup_{r\rightarrow\infty}\frac{log^{+}log^{+}T(r,f)}{logr}$$
by the order  and the hyper-order  of $f$, respectively.

Let $f(z)$ be a meromorphic function, and a finite complex number $\eta$, we define  its difference operators by
\begin{equation*}
\Delta_\eta f(z)=f(z+\eta)-f(z), \quad \Delta_\eta^{n}f(z)=\Delta_{\eta}^{n-1}(\Delta_\eta f(z)).
\end{equation*}

A meromorphic function $a$  satisfying $T(r,a)=S(r,f)$ is called a small function of $f$. We say that two nonconstant meromorphic functions  $f$ and $g$ share small function $a$ CM(IM) if $f-a$ and $g-a$ have the same zeros counting multiplicities (ignoring multiplicities). And we that $f(z)$ and $g(z)$ share $a$ CM almost if
$$N(r,\frac{1}{f-a})+N(r,\frac{1}{g-a})-2N(r,f=a=g)=S(r,f)+S(r,g).$$

In 1977, Rubel and Yang \cite{ruy}  considered the uniqueness of an entire function and its derivative. They proved.

\

{\bf Theorem A} \ Let $f$ be a non-constant entire function, and let $a, b$ be two finite distinct complex values. If $f(z)$ and $f$
 share $a, b$ CM, then $f\equiv f'$.

In recent years, there has been tremendous interests in developing  the value distribution of meromorphic functions with respect to difference analogue,  see [2-8, 10-14, 19]. Heittokangas et al \cite{hkl} proved a similar result analogue of Theorem A concerning shift.

\

{\bf Theorem B}
 Let $f(z)$ be a non-constant entire function of finite order, let $\eta$ be a nonzero finite complex value, and let $a, b$ be two finite distinct complex values.
If $f(z)$ and $f(z+\eta)$ share $a, b$ CM, then $f(z)\equiv f(z+\eta).$

Recently,  Chen-Yi  \cite{cy}, Zhang-Liao \cite{zl},  and Liu-Yang-Fang \cite{lyf}  proved

\

{\bf Theorem C}
 Let $f$ be a transcendental entire function of finite order,  let $\eta$ be a non-zero complex number, $n$ be a positive integer,  and let $a, b$ be two distinct small functions of $f$. If $ f$ and $\Delta_{\eta}^{n}f$ share $a$, $b$ CM, then $ f\equiv \Delta_{\eta}^{n}f$.

In 2014, Halburd-Korhonen-Tohge \cite{h3} investigated the relationship of characteristic functions between $f(z)$ and $f(z+\eta)$ in $\rho_{2}(f)<1$. They obtain the following Lemma 2.1. Immediately, Theorem B and Theorem C are still true when finite order is replaced by $\rho_{2}(f)<1$.

 Li-Yi-Kang \cite{lyk}, L$\ddot{u}$-L$\ddot{u}$ \cite{ll}, Gao, et.al \cite{gkzz} improved Theorem C to meromorphic function. They proved.

\

{\bf Theorem D}
Let $f$ be a  transcendental meromorphic function of hyper-order less than $1$, let $\eta$ be a non-zero complex number, $n\geq1$, an integer, and let $a,b,c$  be three  distinct periodic small functions $f$ with period $\eta$. If $f$ and $\Delta_{\eta}^{n}f$ share $a,b,c$ CM, then $f\equiv\Delta_{\eta}^{n}f$.

In this paper, we prove.

\

{\bf Theorem 1}
Let $f$ be a  non-constant meromorphic function of hyper-order less than $1$, let $\eta$ be a non-zero complex number, $n\geq1$, an integer, and let $a,b,c\in\hat{S}(f)$  be three  distinct small functions and two of them  be periodic small functions with period $\eta$. If $f$ and $\Delta_{\eta}^{n}f$ share $a,b,c$ CM, then $f\equiv\Delta_{\eta}^{n}f$.

\section{Some Lemmas}
\begin{lemma}\label{21l}\cite{h1} Let $f$ be a non-constant meromorphic function of $\rho_{2}(f)<1$,  and let $\eta$ be a non-zero complex number. Then
$$m(r,\frac{f(z+\eta)}{f(z)})=S(r, f),$$
for all r outside of a possible exceptional set E with finite logarithmic measure.
\end{lemma}

\begin{lemma}\label{22l}\cite{h1,h2}
Let $f$ be a non-constant meromorphic function of $\rho_{2}(f)<1$, and let $\eta\neq0$ be a finite complex number. Then
$$T(r,f(z+\eta))=T(r, f(z))+S(r,f).$$
\end{lemma}

\begin{lemma}\label{23l}\cite{y3} Let $f(z)$  be a non-constant meromorphic function of $\rho_{2}(f)<1$, and let $a,b,c$, be three distinct small functions of $f$. Then
$$T(r,f)\leq \overline{N}(r,\frac{1}{f-a})+\overline{N}(r,\frac{1}{f-b})+\overline{N}(r,\frac{1}{f-c})+S(r,f).$$
\end{lemma}

\begin{lemma}\label{24l} Let $f(z)$  be a  transcendental meromorphic function of $\rho_{2}(f)<1$, $n$ a positive integer, let $\eta$ be a nonzero finite number, and let $a(z)\not\equiv\infty, b(z)\not\equiv\infty$ and $c(z)$ be three distinct small meromorphic functions of $f(z)$. Suppose
\[L(f)=\left|\begin{array}{rrrr}a-b& &f-a \\
a'-b'& &f'-a'\end{array}\right|\]
then $L(f)\not\equiv0$.
\end{lemma}
\begin{proof}
Suppose that $L(f)\equiv0$, then we can get $\frac{f'-a'}{f-a}\equiv\frac{a'-b'}{a-b}$. Integrating both side of above we can obtain $f-a=C_{1}(a-b)$, where $C_{1}$ is a nonzero constant. So we have $T(r,f)=S(r,f)$, a contradiction. Hence $L(f)\not\equiv0$.
\end{proof}

\begin{lemma}\label{25l}  Let $f(z)$  be a  transcendental meromorphic function of $\rho_{2}(f)<1$, $n$ a positive integer, let $\eta$ be a nonzero finite number, and let $a(z)\not\equiv\infty, b(z)\not\equiv\infty$ and $c(z)$ be three distinct small meromorphic functions of $f(z)$. Then
$$m(r,\frac{L(f)}{f-a})=S(r,f), \quad m(r,\frac{L(f)}{f-b})=S(r,f).$$
And
$$m(r,\frac{L(f)f}{(f-a)(f-b)})=S(r,f),$$
where $L(f)$ is defined as in Lemma 2.4.
\end{lemma}
\begin{proof}
Obviously, we have
$$m(r,\frac{L(f)}{f-a})\leq m(r,\frac{(a'-b')(f-a)}{f-a})+m(r,\frac{(a-b)(f'-a')}{f-a})=S(r,f).$$
As $\frac{L(f)f}{(f-a)(f-b)}=\frac{C_{1}L(f)}{f-a}+\frac{C_{2}L(f)}{f-b}$,
where $C_{i}(i=1,2)$ are small functions of $f$. Thus
$$m(r,\frac{L(f)f}{(f-a)(f-b)})\leq m(r,\frac{C_{1}L(f)}{f-a})+m(r,\frac{C_{2}L(f)}{f-b})=S(r,f).$$
\end{proof}

\begin{lemma}\label{27l}\cite{y1}
Let $f$ and $g$ be two non-constant meromorphic functions. If $f$ and $g$ share $0,1,\infty$ CM, then
$$N_{(2}(r,f)+N_{(2}(r,\frac{1}{f})+N_{(2}(r,\frac{1}{f-1})=S(r,f).$$
\end{lemma}

\begin{lemma}\label{28l}\cite{y1}
 Let $f$ and $g$ be two non-constant meromorphic functions. If $f$ and $g$ share $0,1,\infty$ CM, and $f$ is not a M$\ddot{o}$bius transformation of $g$,  then\\
(i) $T(r,f)=N(r,\frac{1}{g'})+N_{0}(r)+S(r,f), T(r,g)=N(r,\frac{1}{f'})+N_{0}(r)+S(r,f)$, where $N_{0}(r)$ denotes the zeros of $f-g$, but not the zeros of $f$, $f-1$, and $\frac{1}{f}$.\\
(ii) $T(r,f)+T(r,g)=N(r,f)+N(r,\frac{1}{f})+N(r,\frac{1}{f-1})+N_{0}(r)+S(r,f)$;\\
(iii) $T(r,f)=N(r,\frac{1}{f-a})+S(r,f),$ where $a\neq0,1,\infty$.
\end{lemma}

\begin{lemma}\label{29l}\cite{y1}
Let $f$ and $g$ be two non-constant meromorphic functions. If $f$ and $g$ share $0,1,\infty$ CM with finite lower order, then $T(r,f)=T(r,g)+S(r,f)$.
\end{lemma}

\begin{lemma}\label{2010}\cite{y1} Let $f$ and $g$ be two non-constant meromorphic functions. If $f$ and $g$ share $0,1,\infty$ CM, and
$$N(r,f)\neq T(r,f)+S(r,f),$$
$$N(r,\frac{1}{f-a})\neq T(r,f)+S(r,f),$$
where $a\neq0,1,\infty$. Then $a,\infty$ are the Picard exceptional values of $f$, and $1-a,\infty $ are the Picard exceptional values of $g$.
\end{lemma}

\begin{lemma}\label{2011} \cite{a}
 Let $f$ and $g$ be two nonconstant meromorphic functions. If $f$ and $g$ share $0,1,\infty$ CM, and $f$ is  a M$\ddot{o}$bius transformation of $g$,  then $f$ and $g$ assume one of the following six relations: (i) $fg=1$; (ii) $(f-1)(g-1)=1$; (iii) $f+g=1$; (iv) $f=cg$; (v) $f-1=c(g-1)$; (vi) $[(c-1)f+1][(c-1)g-c]=-c$, where $c\neq0,1$ is a complex number.
\end{lemma}

\begin{lemma}\label{2010}
Let $f$ and $g$ be two non-constant meromorphic functions satisfying
$$\overline{N}(r,f)+\overline{N}(r,g)+\overline{N}(r,\frac{1}{f})+\overline{N}(r,\frac{1}{f})=S(r,f).$$
If $f^{s}g^{t}\equiv1$ for all integers $s$ and $t$($|s|+|t|>0$), then for any positive number $\varepsilon$, we have
$$N_{0}(r,1;f;g)\leq\varepsilon (T(r,f)+T(r,g))+S(r),$$
where  $N_{0}(r,1;f;g)$ denotes the reduced counting function of $f$ and $g$ related to the common $1$-points and $S(r)=o(T(r,f)+T(r,g))$ as $r\rightarrow\infty, r\not\in E_{3}$
 \end{lemma}

\section{The proof of Theorem 1 }
Let $g=\Delta_{\eta}^{n}f$.  Suppose $f\not\equiv g$. Without lose of generality, we discuss two cases, i.e. $c\equiv\infty$ and $c\not\equiv\infty$.

{\bf Case 1} $c\equiv\infty$. Since $f$ is a non-constant meromorphic function satisfying $\rho_{2}(f)<1$, and $f$ and $g$ share $a,b,\infty$ CM, we know that there are two entire functions $p_{1}$ and $p_{2}$ such that
\begin{align}
\frac{g-a}{f-a}=e^{p_{1}}, \quad \frac{g-b}{f-b}=e^{p_{2}}.
\end{align}

Set
\begin{align}
\varphi=\frac{L(f)(f-g)}{(f-a)(f-b)},
\end{align}
where $L(f)\not\equiv0$ is defined as in Lemma 2.4. Since $f\not\equiv g$, then $\varphi\not\equiv0$.

 Set $F=\frac{f-a}{b-a}$ and $G=\frac{g-a}{b-a}$, and thus $F$ and $G$ share $0,1,\infty$ CM, as $f$ and $g$ share $a,b,\infty$ CM. Then by Lemma 2.6, we have
\begin{align}
N(r,f)=N_{1}(r,f),  N(r,\frac{1}{f-a})=N_{1}(r,\frac{1}{f-a}), N(r,\frac{1}{f-b})=N_{1}(r,\frac{1}{f-b}).
\end{align}

Since $f$ is a non-constant  meromorphic function satisfying $\rho_{2}(f)<1$, by Lemma 2.8, we have
\begin{align}
T(r,f)=T(r,F)+S(r,f)=T(r,G)+S(r,f)=T(r,g)+S(r,f).
\end{align}
We claim that
\begin{align}
T(r,f)=N(r,f)+S(r,f).
\end{align}
Otherwise, by Lemma 2.9, we know $N(r,f)=S(r,f)$, and hence {\bf Remark 1} implies $f\equiv g$, a contradiction. We also claim that $F$ is not a M$\ddot{o}$bius transformation of $G$. Otherwise, by Lemma 2.10, if  (i) occurs, we can see that
 \begin{align}
 N(r,\frac{1}{f-a})= N(r,\frac{1}{g-a})=S(r,f), N(r,f)= N(r,g)=S(r,f).
\end{align}
Then by Theorem C, we can obtain a contradiction.\\

If (ii) occurs, we can see that
 \begin{align}
 N(r,\frac{1}{f-b})= N(r,\frac{1}{g-b})=S(r,f), N(r,f)= N(r,g)=S(r,f).
\end{align}
Then by  {\bf Remark 1}, we can obtain a contradiction.

If (iii) occurs, then by (3.1), we can get $f+g=a+b$, that is
\begin{align}
 N(r,\frac{1}{f-a})= N(r,\frac{1}{g-a})=S(r,f),  N(r,\frac{1}{f-b})= N(r,\frac{1}{g-b})=S(r,f),
  \end{align}
and with a similar method of proving {\bf Case 2 (i)}, we can obtain a contradiction since one of $a$ and $b$ is a periodic small function.

If (iv) occurs, that is $F=jG$, where $j\neq0, 1$ is a finite constant. And hence $e^{p_{1}}=j$. So by Lemma 2.1 and (3.1), we have
\begin{align}
T(r,f)&=m(r,\frac{1}{f-b})+S(r,f)=T(r,e^{p_{2}})+S(r,f).
\end{align}
It follows from above that
\begin{eqnarray*}
\begin{aligned}
T(r,f)&=m(r,\frac{1}{f-b})+S(r,f)\notag\\
&\leq m(r,\frac{g-\Delta_{\eta}b}{f-b})+m(r,\frac{1}{g-\Delta_{\eta}b})+S(r,f)\notag\\
&\leq T(r,g)-N(r,\frac{1}{g-\Delta_{\eta}b})+S(r,f),
\end{aligned}
\end{eqnarray*}
that is
\begin{align}
N(r,\frac{1}{g-\Delta_{\eta}^{n}b})=S(r,f).
\end{align}
We claim that $\Delta_{\eta}b\not\equiv b$. Otherwise, by Lemma 2.1 and (3.1), we can obtain
\begin{align}
T(r,e^{p_{2}})=m(r,e^{p_{2}})=m(r,\frac{g-\Delta_{\eta}^{n}b}{f-b})=S(r,f).
\end{align}
Then by (3.9) and (3.11), we have $T(r,f)=T(r,F)+S(r,f)=S(r,f)$, a contradiction.\\
Rewrite $F=jG$ as $f-(a+j(b-a))=j(g-b)$ and $f-(a+j(\Delta_{\eta}^{n}b-a))=j(g-\Delta_{\eta}^{n}b)$.  So
\begin{align}
N(r,\frac{1}{f-a-j(b-a)})=S(r,f),\quad N(r,\frac{1}{f-a-j(\Delta_{\eta}^{n}b-a)})=S(r,f).
\end{align}
Since $j\neq0,1$ and $\Delta_{\eta}^{n}b\not\equiv b$, we know $a+j(b-a)\not\equiv b$ and $a+j(b-a)\not\equiv a+j(\Delta_{\eta}^{n}b-a)$. On the other hand, if $a+j(\Delta_{\eta}^{n}b-a)\not\equiv b$.  Then it follows from  Lemma 2.3, and (3.12) we can get
\begin{align}
T(r,f)&\leq N(r,\frac{1}{f-b})+N(r,\frac{1}{f-a-j(b-a)})\notag\\
&+N(r,\frac{1}{f-a-j(\Delta_{\eta}^{n}b-a)})\leq +S(r,f)=S(r,f).
\end{align}
It is impossible. Hence $a+j(\Delta_{\eta}^{n}b-a)\equiv b$. Set $d=a+j(b-a)$, and we define
\begin{align}
E&=(f-d)(\Delta_{\eta}^{n}(d-b))-(g-\Delta_{\eta}^{n}d)(d-b)\notag\\
&=(f-b)(\Delta_{\eta}^{n}(d-b))-(g-\Delta_{\eta}^{n}b)(d-b).
\end{align}
If $E\not\equiv0$, then by (3.5), (3.12), (3.14) and Lemma 2.1, we have
\begin{eqnarray*}
\begin{aligned}
2T(r,f)&=m(r,\frac{1}{f-b})+m(r,\frac{1}{f-d})+S(r,f)\\
&=m(r,\frac{1}{f-b}+\frac{1}{f-d})+S(r,f)\\
&\leq m(r,\frac{E}{f-b}+\frac{E}{f-d})+m(r,\frac{1}{E})+S(r,f)\\
&\leq T(r,(f-d)(\Delta_{\eta}^{n}(d-b))-(g-\Delta_{\eta}d)^{n}(d-b))+S(r,f)\\
&\leq T(r,f)+S(r,f),
\end{aligned}
\end{eqnarray*}
which is $T(r,f)=S(r,f)$, a contradiction. Therefore $E\equiv0$, i.e.
\begin{align}
(f-d)(\Delta_{\eta}^{n}(d-b))\equiv(g-\Delta_{\eta}^{n}d)(d-b).
\end{align}
Easy to see from  (3.12) and (3.15) that
\begin{align}
N(r,\frac{1}{g-\Delta_{\eta}^{n}d})=S(r,f).
\end{align}

If $\Delta_{\eta}^{n}d\equiv b$, since $c=\infty$, then one of $a,b$ is a  periodic functions, we discuss two cases.

{\bf Case 1.1} $b$ is   periodic. Then from the fact that $f-d=j(g-b)$ and also $d-b=(a-b)(1-j)$, we can know from (3.15) that
\begin{align}
b=\Delta_{\eta}^{n}d=(1-j)(a-b),
\end{align}
and
\begin{align}
b=\Delta_{\eta}^{n}d=(1-j)\Delta_{\eta}^{n}a.
\end{align}
$a+j(\Delta_{\eta}^{n}b-a)\equiv b$ can deduce
\begin{align}
b=\Delta_{\eta}^{n}d=(1-j)a,
\end{align}
it follows from (3.18) that $a$ is also periodic and $\Delta_{\eta}^{n}a\equiv a$. So we have $a=b=0$, a contradiction.

{\bf Case 1.2} $a$ is   periodic. Then from the fact that $f-d=j(g-b)$ and also $d-b=(a-b)(1-j)$, we can know from (3.15) that
\begin{align}
b=\Delta_{\eta}^{n}d=j\Delta_{\eta}^{n}b.
\end{align}
and
\begin{align}
j(b-\Delta_{\eta}^{n}b)=(1-j)(a-b).
\end{align}

(3.20) and (3.21) imply
\begin{align}
a=(1+j)b,
\end{align}
and hence we can obtain from (3.20) and (3.22) that $a=b=0$, a contradiction.

If $\Delta_{\eta}^{n}d\equiv \Delta_{\eta}^{n}b$, then we can obtain  from (3.13) and $j\neq 0,1$ that $T(r,f)=T(r,g)+S(r,f)=T(r,\Delta_{\eta}^{n}d)=S(r,f)$, a contradiction. \\
By Lemma 2.3 and (3.4), we have
\begin{eqnarray*}
\begin{aligned}
T(r,f)&=T(r,g)+S(r,f)\leq N(r,\frac{1}{g-b})+N(r,\frac{1}{g-\Delta_{\eta}^{n}b})\\
&+N(r,\frac{1}{g-\Delta_{\eta}^{n}d})+S(r,f)=S(r,f),
\end{aligned}
\end{eqnarray*}
which is $T(r,f)=S(r,f)$, a contradiction.

If (v) occurs, that is $F-1=i(G-1)$, where $i\neq0, 1$ is a finite constant. And with a similar method of proving (iv), we can obtain a contradiction.

If (vi) occurs, $[(k-1)F+1][(k-1)G-k]=-k$, where $k\neq0,1$ is a complex number. We can see that
 \begin{align}
 N(r,F)=N(r,f)= N(r,g)=N(r,G)=S(r,f).
\end{align}
Then by  {\bf Remark 1}, we can obtain a contradiction.\\

Hence, $F$ is not a M$\ddot{o}$bius transformation of $G$. If $ab\equiv0$, and without lose of generality, we set $a\equiv0$.
Easy to see from (3.1), Lemma 2.1 and Lemma 2.5 that
\begin{eqnarray*}
\begin{aligned}
T(r,\varphi)&=m(r,\frac{L(f)(f-g)}{(f-a)(f-b)})+N(r,\varphi)\notag\\
&\leq m(r,\frac{L(f)f}{(f-a)(f-b)})+m(r,1-\frac{g}{f})+N(r,\varphi)\notag\\
&\leq N_{1}(r,f)+S(r,f),
\end{aligned}
\end{eqnarray*}
that is
\begin{align}
T(r,\varphi)\leq N_{1}(r,f)+S(r,f).
\end{align}

We also obtain
\begin{align}
m(r,\frac{\varphi}{f})&\leq m(r,\frac{L(f)f}{(f-a)(f-b)})+m(r,1-\frac{g}{f})=S(r,f).
\end{align}
Then it follows from Lemma 2.7, (3.2)-(3.4), and (3.24)-(3.25) that
\begin{eqnarray*}
\begin{aligned}
m(r,\frac{1}{f})&\leq m(r,\frac{\varphi}{f})+m(r,\frac{1}{\varphi})\\
&\leq T(r,\varphi)-N(r,\frac{1}{\varphi})+S(r,f)\\
&\leq T(r,\varphi)-(N(r,\frac{1}{L(f)})+N_{0}(r,\frac{1}{f-g}))+S(r,f)\\
&\leq N_{1}(r,f)-T(r,f)+S(r,f)=S(r,f),
\end{aligned}
\end{eqnarray*}
which is
\begin{align}
m(r,\frac{1}{f})=S(r,f).
\end{align}
Here, $N_{0}(r,\frac{1}{f-g})=N_{0}(r,\frac{1}{F-G})+S(r,f)$. So
\begin{align}
T(r,f)=N(r,\frac{1}{f})+S(r,f).
\end{align}
Combing Lemma 2.7, (3.2)-(3.4) and (3.27), we can get
\begin{eqnarray*}
\begin{aligned}
&N(r,f)+N(r,\frac{1}{f})+N(r,\frac{1}{f-b})+N_{0}(r)\\
&=T(r,f)+T(r,g)+S(r,f)\\
&=N(r,\frac{1}{f})+N(r,f)+S(r,f),
\end{aligned}
\end{eqnarray*}
that is
\begin{align}
N(r,\frac{1}{f-b})+N_{0}(r)=S(r,f),
\end{align}
and therefore by (3.28), we have
\begin{align}
T(r,e^{p_{1}})&=N(r,\frac{1}{e^{p_{1}}-1})+S(r,f)\notag\\
&\leq N_{0}(r)+N(r,\frac{1}{f-b})=S(r,f)
\end{align}
and
\begin{align}
&T(r,f)=m(r,\frac{1}{f-b})+N(r,\frac{1}{f-b})+S(r,f)\notag\\
&=m(r,\frac{1}{f-b})+S(r,f)\leq m(r,\frac{1}{g-\Delta_{\eta}^{n}b})+S(r,f)\notag\\
&\leq T(r,g)-N(r,\frac{1}{g-\Delta_{\eta}^{n}b})+S(r,f),
\end{align}
which implies
\begin{align}
N(r,\frac{1}{g-\Delta_{\eta}^{n}b})=S(r,f).
\end{align}
If $\Delta_{\eta}^{n}b=0$, then (3.27) deduces $T(r,f)=S(r,f)$, a contradiction. Hence $\Delta_{\eta}^{n}b=b$. Then by (3.1) and Lemma 2.1, we have
\begin{align}
m(r,e^{p_{2}})=m(r,\frac{g-\Delta_{\eta}^{n}b}{f-b})=S(r,f).
\end{align}
Solving the equation (3.1), we can get
\begin{align}
f=\frac{a-b+be^{p_{2}}-ae^{p_{1}}}{e^{p_{2}}-e^{p_{1}}}.
\end{align}
It follows from above, (3.29) and (3.32),  we have $T(r,f)=S(r,f)$, a contradiction. Hence, we know that neither $\Delta_{\eta}b=0$ nor $\Delta_{\eta}b=b$ holds. Then by Lemma 2.1-Lemma 2.3 and (3.21), we have
\begin{align}
N(r,\frac{1}{f-b})&=N(r,\frac{1}{\Delta_{\eta}^{n}f-b})+S(r,f)\notag\\
&=N(r,\frac{1}{\Delta_{\eta}^{n}f-\Delta_{\eta}^{n}b})+S(r,f)=S(r,f),
\end{align}
and
\begin{align}
T(r,f)=N(r,\frac{1}{f})+S(r,f)=N(r,f)+S(r,f).
\end{align}

Set
$$P_{1}=\frac{\Delta_{\eta}^{n}f-b}{f-b},\quad Q_{1}=\frac{(f-b)\Delta_{\eta}^{n}b}{(\Delta_{\eta}^{n}f-\Delta_{\eta}^{n}b)b}.$$
We ca see from above that
$$N(r,P_{1})+N(r,Q_{1})+N(r,\frac{1}{P_{1}})+N(r,\frac{1}{Q_{1}})=S(r,f).$$
If for all integers $s$ and $t$ satisfying ($|s|+|t|>0$) such that $P_{1}^{s}Q_{1}^{t}\equiv1$, then by Lemma 2.12, we get
\begin{align}
T(r,f)&=N(r,\frac{1}{f-a})+S(r,f)\leq\varepsilon(T(r,P_{1})+T(r,Q_{1}))+S(r,f)\notag\\
&\leq2\varepsilon T(r,f)+S(r,f),
\end{align}
it follows from above and $\varepsilon<\frac{1}{2}$ that $T(r,f)=S(r,f)$, a contradiction. Therefore, there exist two integer $s=1$ and $t=1$ such that $P_{1}Q_{1}\equiv1$. That is
$$\frac{(\Delta_{\eta}^{n}f-b)\Delta_{\eta}^{n}b}{(\Delta_{\eta}^{n}f-\Delta_{\eta}^{n}b)b}\equiv1.$$
Rewrite above as
\begin{align}
\frac{\Delta_{\eta}^{n}b-b}{\Delta_{\eta}^{n}f-\Delta_{\eta}^{n}b}\equiv \frac{b}{\Delta_{\eta}^{n}b}-1,
\end{align}
which follows from Lemma 2.1, Lemma 2.2 and (3.21) that
\begin{align}
T(r,f)=m(r,\frac{1}{f-b})+S(r,f)\leq m(r,\frac{1}{\Delta_{\eta}^{n}f-\Delta_{\eta}^{n}b})+S(r,f)=S(r,f),
\end{align}
but it is impossible.

So $ab\not\equiv0$. By Lemma 2.7 and (3.5) that
\begin{eqnarray*}
\begin{aligned}
3T(r,f)+N_{0}(r)=2T(r,f)+m(r,\frac{1}{f-a})+m(r,\frac{1}{f-b})+S(r,f),
\end{aligned}
\end{eqnarray*}
which follows from above inequality that
\begin{align}
T(r,f)+N_{0}(r)&=m(r,\frac{1}{f-a})+m(r,\frac{1}{f-b})+S(r,f)\notag\\
&\leq m(r,\frac{g-\Delta_{\eta}^{n}a}{f-a})+m(r,\frac{g-\Delta_{\eta}^{n}b}{f-b})+m(r,\frac{1}{g-\Delta_{\eta}^{n}a})+m(r,\frac{1}{g-\Delta_{\eta}^{n}b})\notag\\
&+S(r,f)\leq m(r,\frac{1}{g-\Delta_{\eta}^{n}a})+m(r,\frac{1}{g-\Delta_{\eta}^{n}b})+S(r,f).
\end{align}

We discuss two case.

{\bf Case 1.1} $\Delta_{\eta}^{n}a\not\equiv\Delta_{\eta}^{n}b$.

{\bf Case 1.1.1} $\Delta_{\eta}^{n}a\not\equiv a,b$ and $\Delta_{\eta}^{n}b\not\equiv a,b$. We discuss two case.

{\bf Case 1.1} $\Delta_{\eta}^{n}a\not\equiv\Delta_{\eta}^{n}b$.

{\bf Case 1.1.1} $\Delta_{\eta}^{n}a\not\equiv a,b$ and $\Delta_{\eta}^{n}b\not\equiv a,b$. Let $F_{1}=\frac{1}{F}$ and $G_{1}=\frac{1}{G}$. We only need to discuss $F_{1}$ is not a M$\ddot{o}$bius transformation of $G_{1}$.  We discuss two subcases.

{\bf Subcases 1.1} \quad $T(r,F_{1})\neq N(r,F_{1})+S(r,f)=N(r,\frac{1}{f-a})+S(r,f)$.  Then by Lemma 2.8, we have
\begin{align}
N(r,\frac{1}{f-a})=N(r,\frac{1}{\Delta_{\eta}^{n}f-a})=S(r,f).
\end{align}
Lemma 2.1 implies that
\begin{align}
N(r,\frac{1}{\Delta_{\eta}^{n}f-\Delta_{\eta}^{n}a})=S(r,f).
\end{align}
If $a\equiv a_{\eta}$, then by Lemma 2.1 and (3.1), we get
\begin{align}
m(r,e^{p})=m(r,\frac{\Delta_{\eta}^{n}f-\Delta_{\eta}^{n}a}{f-a})=S(r,f),
\end{align}
and then
\begin{align}
N(r,\frac{1}{f-b})+N_{0}(r)\leq N(r,\frac{1}{e^{p}-1})+S(r,f)=S(r,f).
\end{align}
(3.31) and (3.33) deduce
\begin{align}
m(r,e^{q})&=N(r,\frac{1}{e^{q}-1})+S(r,f)\leq N(r,\frac{1}{f-a})+N_{0}(r)\notag\\
&+S(r,f)=S(r,f).
\end{align}
Combining (3.23), (3.32) and (3.34), we have $T(r,f)=S(r,f)$, a contradiction. That is to say $\Delta_{\eta}^{n}a\not\equiv\infty$.

Set
$$P_{2}=\frac{\Delta_{\eta}^{n}f-a}{f-a},\quad Q_{2}=\frac{(f-a)(b-\Delta_{\eta}^{n}a)}{(\Delta_{\eta}^{n}f-\Delta_{\eta}^{n}a)(b-a)}.$$
We ca see from above that
$$N(r,P_{1})+N(r,Q_{1})+N(r,\frac{1}{P_{1}})+N(r,\frac{1}{Q_{1}})=S(r,f).$$
With a similar method, we can obtain $T(r,f)=S(r,f)$, a contradiction.

{\bf Subcases 1.2} \quad $T(r,F_{1})= N(r,F_{1})+S(r,f)=N(r,\frac{1}{f-a})+S(r,f)$. It follows from (3.29) that
\begin{align}
N(r,\frac{1}{f-b})+N_{0}(r)=S(r,f).
\end{align}
Moreover, we have (3.24). As we set $P_{1}$ and $Q_{1}$ of above, we can also obtain $T(r,f)=S(r,f)$, a contradiction.

{\bf Case 1.1.2} $\Delta_{\eta}^{n}a\equiv a$ and $\Delta_{\eta}^{n}b\equiv b$. Then by Lemma 2.1 and (3.34), we can get $T(r,f)=T(r,F)+S(r,f)=S(r,f)$, a contradiction.

{\bf Case 1.1.3} $\Delta_{\eta}^{n}a\equiv b$ and $\Delta_{\eta}^{n}b\equiv a$.  We can see from the fact one of $a$ and $b$ is a periodic small function that $a\equiv b=0$, a contradiction.

{\bf Case 1.2} $\Delta_{\eta}^{n}a\equiv\Delta_{\eta}^{n}b$. Then by (3.34) we have
\begin{align}
T(r,f)+N_{0}(r)&=m(r,\frac{1}{f-a})+m(r,\frac{1}{f-b})+S(r,f)\notag\\
&\leq m(r,\frac{1}{f-a}+\frac{1}{f-b})+S(r,f)\notag\\
&\leq m(r,\frac{g-\Delta_{\eta}^{n}a}{f-a}+\frac{g-\Delta_{\eta}^{n}b}{f-b})+m(r,\frac{1}{g-\Delta_{\eta}^{n}a})\notag\\
&+S(r,f)\leq m(r,\frac{1}{g-\Delta_{\eta}^{n}a})+S(r,f)\notag\\
&\leq T(r,g)-N(r,\frac{1}{g-\Delta_{\eta}^{n}a})+S(r,f),
\end{align}
it deduces that
\begin{align}
N(r,\frac{1}{g-\Delta_{\eta}^{n}a})+N_{0}(r)=S(r,f).
\end{align}
It follows from Lemma 2.7 that $\Delta_{\eta}^{n}a\equiv a$ or $\Delta_{\eta}^{n}a=b$.  If  $\Delta_{\eta}^{n}a\equiv a$, then by Lemma 2.1 and (3.1) that
\begin{align}
T(r,e^{p_{1}})=m(r,e^{p_{1}})=m(r,\frac{g-\Delta_{\eta}^{n}a}{f-a})=S(r,f).
\end{align}
On the other hand, by Nevanlinna's Second Fundamental Theorem and (3.36), we have
\begin{align}
T(r,e^{p_{2}})&\leq N(r,\frac{1}{e^{p_{2}}-1})+S(r,f)\notag\\
&\leq N(r,\frac{1}{f-a})+N_{0}(r)=S(r,f).
\end{align}
By (3.33), (3.37) and (3.38) that  that $T(r,f)=S(r,f)$, a contradiction.\\

If  $\Delta_{\eta}^{n}b\equiv \Delta_{\eta}^{n}a\equiv b$, then using a similar proof of above, we can also obtain a contradiction.

{\bf Case 2}  $c\not\equiv\infty$.  Without loss of generality, we suppose $a$ and $c$ to be two periodic small functions with period $\eta$.  Since $f$ and $g$ share $a,b,c$ CM, we set
$$F_{2}=\frac{f-a}{f-b}\cdot\frac{c-b}{c-a}, G_{2}=\frac{g-a}{g-b}\cdot\frac{c-b}{c-a},$$
and
\begin{align}
\frac{(f-b)(g-a)}{(f-a)(g-b)}=e^{h_{1}},  \frac{(f-c)(g-a)}{(f-a)(g-c)}=e^{h_{2}}, \frac{(f-b)(g-c)}{(f-c)(g-b)}=e^{h_{3}}.
\end{align}
And we know that $F_{2}$ and $G_{2}$ share $0,1,\infty$ CM almost. And we also have
\begin{align}
 N(r,\frac{1}{f-a})&=N_{1}(r,\frac{1}{f-a}),\notag\\
 N(r,\frac{1}{f-b})&=N_{1}(r,\frac{1}{f-b}),\notag\\
 N(r,\frac{1}{f-c})&=N_{1}(r,\frac{1}{f-c}).
\end{align}
We  claim that $F_{1}$ is not a M$\ddot{o}$bius transformation of $G_{1}$. Otherwise, then by Lemma 2.10, if (i) occurs, we can see that
 \begin{align}
 N(r,\frac{1}{f-a})= N(r,\frac{1}{g-a})=S(r,f), N(r,\frac{1}{f-b})= N(r,\frac{1}{g-b})=S(r,f).
\end{align}
Then by Lemma 2.1 and (3.3), we have
\begin{align}
 2T(r,f)&=m(r,\frac{1}{f-a})+m(r,\frac{1}{f-b})+S(r,f)\notag\\
 &=m(r,\frac{g-\Delta_{\eta}^{n}a}{f-a})+m(r,\frac{g-\Delta_{\eta}^{n}b}{f-b})+m(r,\frac{1}{g-\Delta_{\eta}^{n}a})\notag\\
 &+m(r,\frac{1}{g-\Delta_{\eta}^{n}b})+S(r,f)\notag\\
& \leq m(r,\frac{1}{g-\Delta_{\eta}^{n}a})+m(r,\frac{1}{g-\Delta_{\eta}^{n}b})+S(r,f)\leq 2T(r,g)\notag\\
 &-N(r,\frac{1}{g-\Delta_{\eta}^{n}a})-N(r,\frac{1}{g-\Delta_{\eta}^{n}b})+S(r,f),
\end{align}
which implies
 \begin{align}
 N(r,\frac{1}{g-\Delta_{\eta}^{n}a})+N(r,\frac{1}{g-\Delta_{\eta}^{n}b})=S(r,f).
  \end{align}
Then we can know from Lemma 2.3, (3.41) and (3.43) that $\Delta_{\eta}^{n}a\equiv a$ or $\Delta_{\eta}^{n}a\equiv b$ and $\Delta_{\eta}^{n}b\equiv a$ or $\Delta_{\eta}^{n}b\equiv b$. If one of $\Delta_{\eta}^{n}a\equiv \Delta_{\eta}^{n}b\equiv a$ and $\Delta_{\eta}^{n}a\equiv \Delta_{\eta}^{n}b\equiv b$ occurs, then by Lemma 2.1 and (3.41), we know that
\begin{align}
 2T(r,f)&=m(r,\frac{1}{f-a})+m(r,\frac{1}{f-b})+S(r,f)\notag\\
 &=m(r,\frac{1}{f-a}+\frac{1}{f-b})+S(r,f)\notag\\
  &=m(r,\frac{g-\Delta_{\eta}^{n}a}{f-a}+\frac{g-\Delta_{\eta}^{n}a}{f-b})+m(r,\frac{1}{g-\Delta_{\eta}^{n}a})\notag\\
 & \leq m(r,\frac{1}{g-\Delta_{\eta}^{n}a})+S(r,f)\leq T(r,g)+S(r,f),
\end{align}
which implies $T(r,f)=S(r,f)$, a contradiction. Hence $\Delta_{\eta}^{n}a\not\equiv \Delta_{\eta}^{n}b$. If $\Delta_{\eta}^{n}a\equiv a$, and $\Delta_{\eta}^{n}b\equiv b$, we set
\begin{align}
D&=(f-a)(\Delta_{\eta}^{n}a-\Delta_{\eta}^{n}b)-(g-\Delta_{\eta}^{n}a)(a-b)\notag\\
&=(f-b)(\Delta_{\eta}^{n}a-\Delta_{\eta}^{n}b)-(g-\Delta_{\eta}^{n}b)(a-b).
\end{align}
We claim that $D\not\equiv0$. Otherwise,  by the equalities $\Delta_{\eta}^{n}a\equiv a$, $\Delta_{\eta}^{n}b\equiv b$, and the definition of $D$, we can get $f\equiv g$, a contradiction. So $D\not\equiv0$. Then it is easy to see that
\begin{align}
 2T(r,f)&=m(r,\frac{1}{f-a})+m(r,\frac{1}{f-b})+S(r,f)\notag\\
 &=m(r,\frac{1}{f-a}+\frac{1}{f-b})+S(r,f)\notag\\
  &=m(r,\frac{D}{f-a}+\frac{D}{f-b})+m(r,\frac{1}{D})\notag\\
 & \leq T(r,f-g)-N(r,\frac{1}{f-g})+S(r,f)\notag\\
 &\leq2T(r,f)-N(r,\frac{1}{f-g})+S(r,f),
 \end{align}
which implies
\begin{align}
 N(r,\frac{1}{f-g})=S(r,f).
\end{align}
Then by Lemma 2.3 and the fact that $f$ and $g$ share $a,b,c$ CM, we have
\begin{align}
T(r,f)&\leq \overline{N}(r,\frac{1}{f-a})+\overline{N}(r,\frac{1}{f-b})+\overline{N}(r,\frac{1}{f-c})+S(r,f)\notag\\
&\leq N(r,\frac{1}{f-g})+S(r,f)=S(r,f),
\end{align}
a contradiction. Thus we know that it must occur that $\Delta_{\eta}^{n}a\equiv b$ and $\Delta_{\eta}^{n}b\equiv a$. We can see from the fact   $a,$ is a  periodic small function that $a\equiv b=0$, a contradiction.

If (ii) occurs, we can see that
 \begin{align}
 N(r,\frac{1}{f-b})= N(r,\frac{1}{g-b})=S(r,f), N(r,\frac{1}{f-c})= N(r,\frac{1}{g-c})=S(r,f).
\end{align}
And similar to the proof of (i), we can obtain a contradiction.

If (iii) occurs, we can see that
 \begin{align}
 N(r,\frac{1}{f-a})= N(r,\frac{1}{g-a})=S(r,f), N(r,\frac{1}{f-c})= N(r,\frac{1}{g-c})=S(r,f).
\end{align}
And similar to the proof of (i), we can obtain a contradiction.

If (iv) occurs, that is $F_{1}=dG_{1}$, and $F_{1}-d=d(G_{1}-1)$, i.e.
 \begin{align}
\frac{f-a}{f-b}=d\frac{g-a}{g-b},
\end{align}
where $d\neq0, 1$ is a finite constant.  It follows from above that
 \begin{align}
N(r,\frac{1}{f-c})=N(r,\frac{1}{g-c})=S(r,f), N(r,\frac{1}{f-h})=S(r,f),
\end{align}
where $h=\frac{(1-d)ab+c(bd-a)}{(1-d)c+da-b}$. If $(1-d)c+da-b\equiv0$, we can obtain that $N(r,f)=S(r,f)$, and we can get a contradiction from {\bf Remark 1}. So with a similar proof of proving (vi), and   $c$ is a periodic small function that  we can obtain a contradiction.

If (v) occurs, we can obtain a contradiction with a similar  proof of (vi).

If (vi) occurs, we can see that
 \begin{align}
 N(r,\frac{1}{f-b})= N(r,\frac{1}{g-b})=S(r,f), N(r,\frac{1}{f-\frac{a(d-1)q+b}{(d-1)q+1}})=S(r,f),
\end{align}
where $q=\frac{c-b}{c-a}$, and $b\not\equiv \frac{a(d-1)q+b}{(d-1)q+1}$. Suppose $s=\frac{a(d-1)q+b}{(d-1)q+1}$, then we can obtain from above and Lemma 2.1 that
\begin{align}
 2T(r,f)&=m(r,\frac{1}{f-b})+m(r,\frac{1}{f-s})+S(r,f)\notag\\
 &=m(r,\frac{g-\Delta_{\eta}^{n}b}{f-b})+m(r,\frac{g-\Delta_{\eta}^{n}s}{f-s})+m(r,\frac{1}{g-\Delta_{\eta}^{n}b})\notag\\
 &+m(r,\frac{1}{g-\Delta_{\eta}^{n}s})+S(r,f)\notag\\
& \leq m(r,\frac{1}{g-\Delta_{\eta}^{n}b})+m(r,\frac{1}{g-\Delta_{\eta}^{n}s})+S(r,f)\leq 2T(r,g)\notag\\
 &-N(r,\frac{1}{g-\Delta_{\eta}^{n}b})-N(r,\frac{1}{g-\Delta_{\eta}^{n}s})+S(r,f),
\end{align}
which implies
 \begin{align}
 N(r,\frac{1}{g-\Delta_{\eta}^{n}b})+N(r,\frac{1}{g-\Delta_{\eta}^{n}s})=S(r,f).
  \end{align}
Then we can know from Lemma 2.3, (3.53) and (3.55) that  $\Delta_{\eta}^{n}b\equiv b$ or $\Delta_{\eta}^{n}b\equiv s$, and $\Delta_{\eta}^{n}s\equiv s$ or $\Delta_{\eta}^{n}b\equiv s$.  If $\Delta_{\eta}^{n}b\equiv \Delta_{\eta}^{n}s$, then
 \begin{align}
  2T(r,f)&=m(r,\frac{1}{f-b})+m(r,\frac{1}{f-s})+S(r,f)\notag\\
 &=m(r,\frac{1}{f-b}+\frac{1}{f-s})+S(r,f)\notag\\
 &=m(r,\frac{g-\Delta_{\eta}^{n}b}{f-b}+\frac{g-\Delta_{\eta}^{n}b}{f-s})+m(r,\frac{1}{g-\Delta_{\eta}^{n}b})+S(r,f)\notag\\
  &\leq m(r,\frac{1}{g-\Delta_{\eta}^{n}b})+S(r,f) \leq T(r,g)+S(r,f),
\end{align}
which implies $T(r,f)=S(r,f)$, a contradiction.  Therefore $\Delta_{\eta}^{n}b\not\equiv\Delta_{\eta}^{n}s$. If $\Delta_{\eta}^{n}b\equiv b$ and $\Delta_{\eta}^{n}s\equiv s$, we set
\begin{align}
H&=(f-b)(\Delta_{\eta}^{n}b-\Delta_{\eta}^{n}s)-(g-b)(b-s)\notag\\
&=(f-s)(\Delta_{\eta}^{n}b-\Delta_{\eta}^{n}s)-(g-s)(b-s).
\end{align}
If $H\equiv0$, we have $f\equiv g$, a contradiction. Hence $H\not\equiv0$.  It is easy to see from (3.67) that
 \begin{align}
  2T(r,f)&=m(r,\frac{1}{f-b})+m(r,\frac{1}{f-s})+S(r,f)\notag\\
 &\leq 2 m(r,\frac{1}{H})+S(r,f) \leq 2T(r,f)-2N(r,\frac{1}{H})+S(r,f),
\end{align}
that is
\begin{align}
N(r,\frac{1}{H})=S(r,f).
\end{align}
Then by Lemma 2.3 and the fact that $f$ and $g$ share $a,b,c$ CM, we have
\begin{align}
T(r,f)&\leq \overline{N}(r,\frac{1}{f-a})+\overline{N}(r,\frac{1}{f-b})+\overline{N}(r,\frac{1}{f-c})+S(r,f)\notag\\
&\leq N(r,\frac{1}{f-g})+S(r,f)=S(r,f),
\end{align}
a contradiction. Hence $\Delta_{\eta}^{n}b\equiv s$ and $\Delta_{\eta}^{n}s\equiv b$. On the other hand, Lemma 2.10 (vi) tells us that $N(r,\frac{1}{g-p})=S(r,f)$, where $p=\frac{aq(d-1)-db}{q(d-1)-d}$. Obviously, $b\not\equiv p$, and $p\equiv s$, otherwise, by Lemma 2.3, we can obtain $a=b$ and $T(r,f)=S(r,f)$, a contradiction.  But $p\equiv s$ implies $d=-1$, and we take it into Lemma 2.10 (vi), we have
$$(2F_{1}-1)(2G_{1}-1)=1,$$
and it follows from above that
$$N(r,\frac{1}{f-(2a-b)})=S(r,f), N(r,\frac{1}{g-(2a-b)})=S(r,f).$$
And then we have
\begin{align}
  2T(r,f)&=m(r,\frac{1}{f-b})+m(r,\frac{1}{f-(2a-b)})+S(r,f)\notag\\
 &=m(r,\frac{1}{f-b}+\frac{1}{f-(2a-b)})+S(r,f)\notag\\
 &=m(r,\frac{g-\Delta_{\eta}^{n}b}{f-b}+\frac{g-\Delta_{\eta}^{n}b}{f-(2a-b)})+m(r,\frac{1}{g-\Delta_{\eta}^{n}b})+S(r,f)\notag\\
  &\leq m(r,\frac{1}{g-\Delta_{\eta}^{n}b})+S(r,f) \leq T(r,g)+S(r,f),
\end{align}
which implies $T(r,f)=S(r,f)$, a contradiction.

Therefore, $F_{1}$ is not a M$\ddot{o}$bius transformation of $G_{1}$.  We discuss two subcases.

 {\bf Subcase 2.1} $T(r,f)\neq N(r,\frac{1}{f-b})+S(r,f)$.  Then by Lemma 2.9 we know that
 \begin{align}
N(r,\frac{1}{f-b})=S(r,f).
\end{align}
Then by Lemma 2.7 and above, we have
\begin{eqnarray*}
\begin{aligned}
3T(r,f)+N_{0}(r)&=2T(r,f)+m(r,\frac{1}{f-a})+m(r,\frac{1}{f-b})+m(r,\frac{1}{f-c})+S(r,f)\\
&= 3T(r,f)+m(r,\frac{1}{f-a})+m(r,\frac{1}{f-c})+S(r,f),
\end{aligned}
\end{eqnarray*}
which implies
 \begin{align}
N_{0}(r)=m(r,\frac{1}{f-a})+m(r,\frac{1}{f-c})+S(r,f).
\end{align}
Furthermore, by Lemma 2.1, we know
\begin{eqnarray*}
\begin{aligned}
T(r,f)&=m(r,\frac{1}{f-b})+S(r,f)\leq m(r,\frac{1}{g-\Delta_{\eta}^{n}b})+S(r,f)\\
&\leq T(r,g)-N(r,\frac{1}{g-\Delta_{\eta}^{n}b})+S(r,f),
\end{aligned}
\end{eqnarray*}
it follows from Lemma 2.8 that
\begin{align}
N(r,\frac{1}{g-\Delta_{\eta}^{n}b})=S(r,f).
\end{align}
And by Lemma 2.8 and Lemma 2.9,
\begin{align}
T(r,f)=N(r,f)+S(r,f)=T(r,g)+S(r,f)=N(r,g)+S(r,f).
\end{align}
We can know from (3.63)-(3.64) and Lemma 2.7 that either $\Delta_{\eta}^{n}b\equiv a$ or $\Delta_{\eta}^{n}b\equiv b$ or $\Delta_{\eta}^{n}b\equiv c$. If $\Delta_{\eta}^{n}b\equiv b$, then Lemma 2.1, Lemma 2.11 and (3.39) deduce
\begin{eqnarray*}
\begin{aligned}
T(r,e^{h_{1}})&=m(r,e^{h_{1}})=m(r,e^{-h_{1}})=m(r,\frac{(f-a)(g-b)}{(f-b)(g-a)})\\
&\leq m(r,\frac{g-b}{f-b})+m(r,\frac{f-a}{g-a})\\
&\leq m(r,\frac{1}{f-a})+N(r,\frac{g-a}{f-a})-N(r,\frac{f-a}{g-a})\\
&\leq  m(r,\frac{1}{f-a})+N(r,g)+N(r,\frac{1}{f-a})-N(r,f)\\
&-N(r,\frac{1}{g-a})+S(r,f)\leq m(r,\frac{1}{f-a})+S(r,f),
\end{aligned}
\end{eqnarray*}
which implies
\begin{align}
T(r,e^{h_{1}})\leq m(r,\frac{1}{f-a})+S(r,f).
\end{align}
We also have
\begin{align}
T(r,e^{h_{3}})\leq m(r,\frac{1}{f-c})+S(r,f).
\end{align}

Applying Lemma 2.3 to $e^{h_{1}}$ and $e^{h_{3}}$, we have
\begin{align}
T(r,e^{h_{1}})=\overline{N}(r,\frac{1}{e^{h_{1}}-1})+S(r,f),
\end{align}
and
\begin{align}
T(r,e^{h_{3}})=\overline{N}(r,\frac{1}{e^{h_{3}}-1})+S(r,f).
\end{align}

It follows from (3.39)-(3.40) and  (3.65)-(3.69) that
\begin{eqnarray*}
\begin{aligned}
\overline{N}(r,\frac{1}{f-c})+\overline{N}_{0}(r)&=\overline{N}(r,\frac{1}{e^{h_{1}}-1})+S(r,f)\\
&\leq m(r,\frac{1}{f-a})+S(r,f),
\end{aligned}
\end{eqnarray*}
which is
\begin{align}
\overline{N}(r,\frac{1}{f-c})+\overline{N}_{0}(r)\leq m(r,\frac{1}{f-a})+S(r,f).
\end{align}
Similarly, we have
\begin{align}
\overline{N}(r,\frac{1}{f-a})+\overline{N}_{0}(r)\leq m(r,\frac{1}{f-c})+S(r,f).
\end{align}
From the fact that the zero of $f-g$ with multiplicity at least $2$ are the zeros of $e^{h_{1}}-1$ with multiplicity at least $2$, and hence we have
\begin{eqnarray*}
\begin{aligned}
2\overline{N}_{0}(r)+\overline{N}(r,\frac{1}{f-a})+\overline{N}(r,\frac{1}{f-c})=N_{0}(r)+S(r,f),
\end{aligned}
\end{eqnarray*}
that is
\begin{align}
\overline{N}_{0}(r)+\overline{N}(r,\frac{1}{f-a})+\overline{N}(r,\frac{1}{f-c})=S(r,f).
\end{align}
Then by Lemma 2.3,  (3.71) and above, we can obtain $T(r,f)=S(r,f)$, a contradiction. Therefore $\Delta_{\eta}^{n}b\equiv a$ or $\Delta_{\eta}^{n}b\equiv c$. If $\Delta_{\eta}^{n}b\equiv a$, then by (3.69) and the fact that $f$ and $g$ share $a$ CM,  we have
\begin{align}
N(r,\frac{1}{f-a})=N(r,\frac{1}{g-a})=S(r,f).
\end{align}
Then similar to the proof of the {\bf Case 2-(i)}, we can obtain a contradiction. If $\Delta_{\eta}^{n}b\equiv c$, we can also get a contradiction.

 {\bf Subcase 2.2} $T(r,f)= N(r,\frac{1}{f-b})+S(r,f)$.  Then by Lemma 2.7 and Lemma 2.8 we know that
 \begin{eqnarray*}
\begin{aligned}
2T(r,f)&=N(r,\frac{1}{f-a})+N(r,\frac{1}{f-b})+N(r,\frac{1}{f-c})+N_{0}(r)+S(r,f)\\
&=T(r,f)+N(r,\frac{1}{f-a})+N(r,\frac{1}{f-c})+N_{0}(r)+S(r,f)\\
&=N(r,\frac{1}{f-g})+S(r,f),
\end{aligned}
\end{eqnarray*}
it follows that
\begin{align}
&T(r,f)=N(r,\frac{1}{f-a})+N(r,\frac{1}{f-c})+N_{0}(r)+S(r,f),\notag\\
&2T(r,f)=N(r,\frac{1}{f-g})+S(r,f).
\end{align}
 That is
 \begin{align}
T(r,f)+N_{0}(r)&=m(r,\frac{1}{f-a})+m(r,\frac{1}{f-c})+S(r,f)\notag\\
&\leq m(r,\frac{g-\Delta_{\eta}^{n}a}{f-a})+m(r,\frac{g-\Delta_{\eta}^{n}c}{f-c})+m(r,\frac{1}{g-\Delta_{\eta}^{n}a})+m(r,\frac{1}{g-\Delta_{\eta}^{n}c})\notag\\
&+S(r,f)\leq m(r,\frac{1}{g-\Delta_{\eta}a})+m(r,\frac{1}{g-\Delta_{\eta}c})+S(r,f).
\end{align}

We discuss two subcase.

{\bf Subcase 2.2.1} $\Delta_{\eta}^{n}a\not\equiv\Delta_{\eta}^{n}c$.

{\bf Subase 2.2.1.1} $\Delta_{\eta}^{n}a\not\equiv a,c$ and $\Delta_{\eta}^{n}c\not\equiv a,c$. Set $F_{3}=\frac{1}{F_{2}}$ and $G_{3}=\frac{1}{G_{2}}$. With the same way to prove {\bf Subcase 2.1}, we only need to discuss $T(r,f)= N(r,\frac{1}{f-a})+S(r,f)$. It follows from (3.79) that
$$N(r,\frac{1}{f-c})+N_{0}(r)=S(r,f).$$
But in this case, we can also obtain a contradiction.

{\bf Subase 2.2.1.2} $\Delta_{\eta}^{n}a\equiv a$ and $\Delta_{\eta}^{n}c\equiv c$. We can get $a\equiv c=0$, a contradiction.

{\bf Subase 2.2.1.3} $\Delta_{\eta}^{n}a\equiv c$ and $\Delta_{\eta}^{n}c\equiv a$.  We can get $a\equiv c=0$, a contradiction.

{\bf Case 2.2.2} $\Delta_{\eta}^{n}a\equiv\Delta_{\eta}^{n}c$. Then by (3.74) we have
\begin{align}
T(r,f)+N_{0}(r)&=m(r,\frac{1}{f-a})+m(r,\frac{1}{f-c})+S(r,f)\notag\\
&\leq m(r,\frac{1}{f-a}+\frac{1}{f-c})+S(r,f)\notag\\
&\leq m(r,\frac{g-\Delta_{\eta}^{n}a}{f-a}+\frac{g-\Delta_{\eta}^{n}c}{f-c})+m(r,\frac{1}{g-\Delta_{\eta}^{n}a})\notag\\
&+S(r,f)\leq m(r,\frac{1}{g-\Delta_{\eta}^{n}a})+S(r,f)\notag\\
&\leq T(r,g)-N(r,\frac{1}{g-\Delta_{\eta}^{n}a})+S(r,f),
\end{align}
it deduces that
\begin{align}
N(r,\frac{1}{g-\Delta_{\eta}^{n}a})+N_{0}(r)=S(r,f).
\end{align}
It follows from Lemma 2.7 that $\Delta_{\eta}^{n}a\equiv a$ or $\Delta_{\eta}^{n}a=c$.  If  $\Delta_{\eta}a^{n}\equiv a$, then with the same proof of (3.75), and by (3.39) and $T(r,f)= N(r,\frac{1}{f-b})+S(r,f)$, we can have
 \begin{align}
T(r,e^{h_{2}})\leq m(r,\frac{1}{f-c})+S(r,f)=S(r,f).
\end{align}
 And thus,
 \begin{align}
\overline{N}_{0}(r)+\overline{N}(r,\frac{1}{f-b})\leq m(r,\frac{1}{f-b})+S(r,f),
\end{align}
which follows from (3.40) that $N(r,\frac{1}{f-b})=S(r,f)$, and furthermore we get $T(r,f)=S(r,f)$, a contradiction. If $\Delta_{\eta}^{n}a=c$, then we can also obtain  a contradiction with a same method of above.

\

{\bf Conflict of Interest}  The author declares that there is  no conflict of interest regarding the publication of this paper.
\

{\bf Acknowledgements} The author would like to thank to anonymous referees for their helpful comments.


\end{document}